\documentclass[11pt,a4paper]{article}%
\usepackage{amsfonts}
\usepackage{amsmath}
\usepackage{amssymb}
\usepackage{graphicx}
\usepackage{geometry}%
\providecommand{\U}[1]{\protect\rule{.1in}{.1in}}
\newtheorem{theorem}{Theorem}

\newtheorem{definition}[theorem]{Definition}

\begin{document}

\title{On Razzaboni Transformation of Surfaces in Minkowski 3-Space }
\author{Melek Erdo\u{g}du\thanks{Department of Mathematics-Computer Sciences,
Necmettin Erbakan University, \ Konya, TURKEY. e-mail: merdogdu@konya.edu.tr.}
\ \thanks{Corresponding author}
\and Mustafa \"{O}zdemir\thanks{ Department of Mathematics, Akdeniz University,
Antalya, TURKEY. e-mail: mozdemir@akdeniz.edu.tr.} }
\maketitle

\begin{abstract}
In this paper, we investigate the surfaces generated by binormal motion of
Bertrand curves, which is called Razzaboni surface, in Minkowski 3-space. We
discussed the geometric properties of these surfaces in $\mathbb{M}^{3}$
according to the character of Bertrand geodesics. Then, we define the
Razzaboni transformation for a given Razzaboni surface. In other words, we
prove that there exists a dual of Razzaboni surface for each case.\ Finally,
we show that Razzaboni transformation maps the surface $\sigma,$ which has
Bertrand geodesics with constant curvature, to the surface $\sigma^{\ast}$
whose Bertrand geodesics also have constant curvature with opposite sign.

\end{abstract}



\section{Introduction}

\noindent Bertrand curves are well-studied classical curves in Euclidean space
and Lorentzian space \cite{eisen}, \cite{lopez}. Curves of constant curvature
or torsion constitute particular Bertrand curves and Bertrand curve shares its
principal normals with another Bertrand curve, which is called Bertrand mate.
The surfaces, which admits geodesic embedding of Bertrand curves, were deeply
studied by Amilcare Razzaboni \cite{raz}. Then the class of surfaces with
Bertrand geodesics came to be called Razzaboni surfaces \cite{main}. In
particular, Razzaboni surfaces, which have Bertrand geodesics with constant
curvature or torsion, was discussed in earlier work of Razzaboni \cite{raze},
\cite{fibbi}.

\bigskip

\noindent The Razzaboni surface can be considered as a surface which is
generated by binormal motion of Bertrand curve \cite{main}. The binormal
motion of curves with constant curvature and torsion is discussed in
\cite{schief}, \cite{backlund}. On the other hand, motion of timelike surfaces
in timelike geodesic coordinates is examined in the study \cite{gurbuz 2}.

\bigskip

\noindent In this study, we investigate the Razzaboni surfaces in Minkowski
3-space in three different cases. The geometric properties of Razzaboni
surfaces in geodesic coordinates are stated for each case. And, we define
Razzaboni transformation for a given Razzaboni surface in $\mathbb{M}^{3}.$
Then, we examine the curvature and torsion of the Bertrand geodesics of dual
Razzaboni surface. Moreover, the Gaussian curvature of Razzaboni surface are
given for each case. Finally, the dual Razzaboni surfaces are examined in the
case geodesic of constant torsion or curvature, respectively.

\section{Preliminaries}

\noindent In this part, we give a brief summary of Razzaboni surface in
Euclidean 3- Space, \cite{main}. Then, some essential properties of curves in
Minkowski 3-space are given to provide the necessary background \cite{I2},
\cite{DI}.

\begin{definition}
Two curves which, at any point, have a common principal normal vector are
called Bertrand curves. Moreover, curves for which there exists constants $A$
and $B$ such that%
\[
A\kappa+B\tau=1
\]
are also known as Bertrand curves. Here $\kappa$ is the curvature and $\tau$
is the torsion of the curve.
\end{definition}

\noindent We know that a curve is a geodesic on a surface $\sigma$ if and only
if the principal normal of the curve is parallel to the normal $N$ of the
surface. This means that if the surface $\sigma$ is spanned by one parameter
family of geodesic Bertrand curves $\alpha$ with the same constants $A$ and
$B$ then the Bertrand mate $\alpha^{\ast}$ form a parallel surface
$\sigma^{\ast}$ on which they are likewise geodesics \cite{main}.

\noindent Let $\sigma=\sigma(u,v)$ be a one parameter family of geodesics with
geodesic coordinates. Then the first fundamental form of the surface is of the
form%
\[
I=du^{2}+\lambda^{2}dv^{2}.
\]
Here, $u-$ parameter curves are unit speed geodesics and $v-$ parameter curves
form the orthogonal parallels. Since $\sigma_{u}\perp\sigma_{v}$ and principal
normal $n$ of the geodesics is orthogonal to the surface, then the tangent
vectors of coordinate line should be of the form%
\[
\sigma_{u}=t,\text{ }\sigma_{v}=\lambda b
\]
where $b$ denotes the binormal vector field of the geodesics. The surface
$\sigma$ is generated by motion of a in extensible curve in binormal direction
with the speed $\lambda.$ Thus, a Razzaboni surface can be considered as a
surface which is generated by binormal motion of Bertrand curve with the same
constants $A$ and $B.$ The variation of Frenet Frame $\{t,n,b\}$ of the
geodesics in $u-$ direction is given by Serret Frenet equations;%
\[
\left[
\begin{array}
[c]{c}%
t\\
n\\
b
\end{array}
\right]  _{u}=\left[
\begin{array}
[c]{ccc}%
0 & \kappa & 0\\
-\kappa & 0 & \tau\\
0 & -\tau & 0
\end{array}
\right]  \left[
\begin{array}
[c]{c}%
t\\
n\\
b
\end{array}
\right]  .
\]
The variation of $\{t,n,b\}$ in $v-$ direction should be of the form
\[
\left[
\begin{array}
[c]{c}%
t\\
n\\
b
\end{array}
\right]  _{v}=\left[
\begin{array}
[c]{ccc}%
0 & \alpha & \beta\\
-\alpha & 0 & \gamma\\
-\beta & -\gamma & 0
\end{array}
\right]  \left[
\begin{array}
[c]{c}%
t\\
n\\
b
\end{array}
\right]  .
\]
By using the compatibility condition $\sigma_{uv}=\sigma_{vu},$ we get
\[
\alpha n+\beta b=-\tau\lambda n+\lambda_{u}b.
\]
Thus, the variation of $\{t,n,b\}$ in $v-$ direction becomes%
\[
\left[
\begin{array}
[c]{c}%
t\\
n\\
b
\end{array}
\right]  _{v}=\left[
\begin{array}
[c]{ccc}%
0 & -\tau\lambda & \lambda_{u}\\
\tau\lambda & 0 & \gamma\\
-\lambda_{u} & -\gamma & 0
\end{array}
\right]  \left[
\begin{array}
[c]{c}%
t\\
n\\
b
\end{array}
\right]  .
\]
The compatibility conditions $t_{uv}=t_{vu},$ $n_{uv}=n_{vu}$ and
$b_{uv}=b_{vu}$ gives the following undetermined system%
\begin{align*}
\kappa_{v}  &  =-2\tau\lambda_{u}-\tau_{u}\lambda,\\
\tau_{v}  &  =\gamma_{u}+\kappa\lambda_{u},\\
\lambda_{uu}  &  =\tau^{2}\lambda+\kappa\gamma.
\end{align*}
The above system can be considered as the Gauss-Minardi-Codazzi equations for
the surface with geodesic coordinates. If the constraint
\[
A\kappa+B\tau=1
\]
is imposed then the system is well determined and the surface is guaranteed to
be a Razzaboni surface \cite{main}.

\noindent In the case of $A=0,$ the geodesics have constant torsion. Choosing
$B=\tau=1$ without loss of generality, the system becomes%
\begin{align*}
\kappa_{v}  &  =-2\lambda_{u},\\
0  &  =\gamma_{u}+\kappa\lambda_{u},\\
\lambda_{uu}  &  =\lambda+\kappa\gamma.
\end{align*}
This integrable system may be regarded as an extension of the sine-Gordon
equation;%
\[
\Psi_{uv}=\sin\Psi.
\]
A single equation is obtained as follows:%
\[
(\frac{\theta_{vuu}-\theta_{v}}{\theta_{u}})_{u}+\theta_{u}\theta_{vu}=0
\]
by setting $\kappa=\theta_{u}$ and $\lambda=-\frac{\theta_{v}}{2}.$

\noindent In the case of $B=0,$ the geodesics have constant curvature.
Choosing $A=\kappa=1$ without loss of generality, the system becomes%
\begin{align*}
0  &  =-2\tau\lambda_{u}-\tau_{u}\lambda,\\
\tau_{v}  &  =\gamma_{u}+\lambda_{u},\\
\lambda_{uu}  &  =\tau^{2}\lambda+\gamma.
\end{align*}
If we set $\lambda=\frac{1}{\sqrt{\tau}},$ then the system reduces to
\[
\tau_{v}=\left[  (\frac{1}{\sqrt{\tau}})_{uu}-\tau\sqrt{\tau}+\frac{1}%
{\sqrt{\tau}}\right]  _{u}%
\]
which represents an extension of the Dym equation%
\[
\tau_{v}=(\frac{1}{\sqrt{\tau}})_{uuu}.
\]
It can be noted that, the extended Dym equation is generated by binormal
motion of an inextensible curve with the speed $\frac{1}{\sqrt{\tau}}.$ For
further discussions about the Razzaboni surfaces, the readers are referred to
\cite{main}.

\noindent The Minkowski 3-space, denoted by $\mathbb{M}^{3},$ is Euclidean
3-space equipped with Lorentzian inner product%
\[
\left\langle u,v\right\rangle _{\mathbb{L}}=-u_{1}v_{1}+u_{2}v_{2}+u_{3}v_{3}%
\]
where $u=(u_{1},u_{2},u_{3}),$ $v=(v_{1},v_{2},v_{3})\in\mathbb{M}^{3}.$
Lorentzian inner product characterizes the elements $u=(u_{1},u_{2},u_{3})$ of
$\mathbb{M}^{3}$.
\[
\text{if }\left\langle u,u\right\rangle _{\mathbb{L}}>0\text{ or }u=0\text{
then }u\text{ is called spacelike},
\]%
\[
\text{if }\left\langle u,u\right\rangle _{\mathbb{L}}<0\text{ then }u\text{ is
called timelike,}%
\]%
\[
\text{if }\left\langle u,u\right\rangle _{\mathbb{L}}=0\text{ and }%
u\neq0\text{ then }u\text{ is called lightlike or null.}%
\]
The norm of $u=(u_{1},u_{2},u_{3})\in\mathbb{M}^{3}$ is defined by $\left\Vert
u\right\Vert =\sqrt{\left\vert \left\langle u,u\right\rangle \right\vert }.$
Lorentzian vector product of $u=(u_{1},u_{2},u_{3})$ and $v=(v_{1},v_{2}%
,v_{3})\in\mathbb{M}^{3}$ is defined by%
\[
u\times_{\mathbb{L}}v=\left\vert
\begin{array}
[c]{ccc}%
-e_{1} & e_{2} & e_{3}\\
u_{1} & u_{2} & u_{3}\\
v_{1} & v_{2} & v_{3}%
\end{array}
\right\vert .
\]
For details, see \cite{I}, \cite{erdo}, \cite{ozdemir}.

\noindent A curve $\alpha$ in $\mathbb{M}^{3}$ is called timelike, spacelike
or null if and only if tangent vector field $t$ of $\alpha$ is timelike,
spacelike or null, respectively. Let $\alpha(s)$ be a unit speed nonlightlike
curve in $\mathbb{M}^{3}$, i,e$.,$ $\left\langle \alpha^{\prime}%
,\alpha^{\prime}\right\rangle _{\mathbb{L}}=\varepsilon_{1}=\pm1.$ The
constant $\varepsilon_{1}$ is called the casual character of $\alpha.$ Frenet
frame field $\{t,n,b\},$ which is an orthonormal field along $\alpha,$ satisfy
the Serret-Frenet equation\textit{:}%
\begin{equation}
\left[
\begin{array}
[c]{c}%
t\\
n\\
b
\end{array}
\right]  _{s}=\left[
\begin{array}
[c]{ccc}%
0 & \varepsilon_{2}\kappa & 0\\
-\varepsilon_{1}\kappa & 0 & -\varepsilon_{3}\tau\\
0 & \varepsilon_{2}\tau & 0
\end{array}
\right]  \left[
\begin{array}
[c]{c}%
t\\
n\\
b
\end{array}
\right]  . \tag{1}%
\end{equation}
The functions $\kappa\geq0$ and $\tau$ are called the curvature and torsion,
respectively. The constants $\varepsilon_{2}=\left\langle n,n\right\rangle
_{\mathbb{L}}$ and $\varepsilon_{3}=\left\langle b,b\right\rangle
_{\mathbb{L}}$ are called the second casual character and third casual
character\textit{ }of $\alpha,$ respectively \cite{I2}, \cite{erdo}, \cite{DI}.

\section{Razzaboni Surfaces Minkowski 3-Space}

\noindent In this part, Razzaboni surfaces in Minkowski 3-space are
investigated in three different cases.

\noindent\textbf{Case 1: Geodesic Bertrand curves of the Razzaboni surface
have timelike principal normal.}

\noindent Let $\sigma=\sigma(u,v)$ be a one parameter family of geodesic
Bertrand curve with timelike principal normals in $\mathbb{M}^{3}.$ And let
$u$ and $v$ be the geodesic coordinates of the Razzaboni surface $\sigma.$
Then the first fundamental form of the surface is of the form%
\[
I=du^{2}+\lambda^{2}dv^{2}.
\]
Here $u-$ parameter curves are unit speed spacelike Bertrand geodesics and
$v-$ parameter curves forms orthogonal spacelike parallels. Since
$\left\langle \sigma_{u},\sigma_{v}\right\rangle _{\mathbb{L}}=0$ and the
principal normal $n$ of geodesics are normal to the surface, then the tangent
vectors of coordinate curves are given as
\[
\sigma_{u}=t\text{ and }\sigma_{v}=\lambda b
\]
where $b$ denotes the binormal vector field of geodesics. The variation of the
Frenet Frame $\{t,n,b,\}$ of Bertrand geodesics in $u-$ direction is obtained
as follows:%
\[
\left[
\begin{array}
[c]{c}%
t\\
n\\
b
\end{array}
\right]  _{u}=\left[
\begin{array}
[c]{ccc}%
0 & -\kappa & 0\\
-\kappa & 0 & -\tau\\
0 & -\tau & 0
\end{array}
\right]  \left[
\begin{array}
[c]{c}%
t\\
n\\
b
\end{array}
\right]
\]
by setting $(\varepsilon_{1},\varepsilon_{2},\varepsilon_{3})=(1,-1,1)$ in
equation 1. The variation of $\{t,n,b\}$ in $v-$ direction must be in the
following form:%
\[
\left[
\begin{array}
[c]{c}%
t\\
n\\
b
\end{array}
\right]  _{v}=\left[
\begin{array}
[c]{ccc}%
0 & \alpha & \beta\\
\alpha & 0 & \gamma\\
-\beta & \gamma & 0
\end{array}
\right]  \left[
\begin{array}
[c]{c}%
t\\
n\\
b
\end{array}
\right]  .
\]
The compatibility condition $\sigma_{uv}=\sigma_{vu}$ implies that%
\[
\alpha n+\beta b=(-\lambda\tau)n+\lambda_{u}b.
\]
So, we get $\alpha=-\lambda\tau$ and $\beta=\lambda_{u}.$Then the variation of
$\{t,n,b\}$ in $v-$ direction becomes
\[
\left[
\begin{array}
[c]{c}%
t\\
n\\
b
\end{array}
\right]  _{v}=\left[
\begin{array}
[c]{ccc}%
0 & -\lambda\tau & \lambda_{u}\\
-\lambda\tau & 0 & \gamma\\
-\lambda_{u} & \gamma & 0
\end{array}
\right]  \left[
\begin{array}
[c]{c}%
t\\
n\\
b
\end{array}
\right]  .
\]
On the other hand, the compatibility condition $t_{uv}=t_{vu}$ implies that
\begin{align*}
t_{uv}  &  =(-\kappa n)_{v}=(\lambda\kappa\tau)t+(-\kappa_{v})n+(-\kappa
\gamma)b,\\
t_{vu}  &  =(-\lambda\tau n+\lambda_{u}b)_{u}=(\lambda\kappa\tau
)t+(-\lambda\tau_{u}-2\lambda_{u}\tau)n+(-\kappa\gamma)b.
\end{align*}
So we get
\[
\kappa_{v}=\lambda\tau_{u}+2\lambda_{u}\tau.
\]
And by the compatibility condition $n_{uv}=n_{vu},$ we obtain%
\begin{align*}
n_{uv}  &  =(-\kappa t-\tau b)_{v}=(-\kappa_{v}+\lambda_{u}\tau)t+(\kappa
\lambda\tau-\tau\gamma)n+(-\lambda_{u}\kappa-\tau_{v})b,\\
n_{vu}  &  =(-\lambda\tau t+\gamma b)_{u}=(-\lambda_{u}\tau-\lambda\tau
_{u})t+(\kappa\lambda\tau-\tau\gamma)n+(\gamma_{u})b.
\end{align*}
Thus, we get
\[
\tau_{v}=-\lambda_{u}\kappa-\gamma_{u}.
\]
Finally, the compatibility condition $b_{uv}=b_{vu}$ gives that%
\begin{align*}
b_{uv}  &  =(-\tau n)_{v}=(\lambda\tau^{2})t+(-\tau_{v})n+(-\tau\gamma)b,\\
b_{vu}  &  =(-\lambda_{u}t+\gamma n)_{u}=(-\lambda_{uu}-\kappa\gamma
)t+(\lambda_{u}\kappa+\gamma_{u})n+(-\tau\gamma)b.
\end{align*}
Then, we get
\[
\lambda_{uu}=-\lambda\tau^{2}-\kappa\gamma
\]
The system%
\begin{align*}
\kappa_{v}  &  =\lambda\tau_{u}+2\lambda_{u}\tau,\\
\tau_{v}  &  =-\lambda_{u}\kappa-\gamma_{u},\\
\lambda_{uu}  &  =-\lambda\tau^{2}-\kappa\gamma
\end{align*}
may be regarded as the Gauss-Minardi-Codazzi equations for the surface. If the
constraint
\[
A\kappa+B\tau=1
\]
is imposed then the system is well determined and the surface $\sigma$ is
guaranteed to be a Razzaboni surface.

\bigskip

\noindent Now, let's suppose that $\varphi$ is a transformation between two
Razzaboni surfaces in $\mathbb{M}^{3}$ defined as $\sigma^{\ast}%
=\varphi\left(  \sigma\right)  .$ And let $\{t^{\ast},n^{\ast},b^{\ast}\}$ be
the Frenet Frame of Bertrand geodesics of Razzaboni surface $\sigma^{\ast}.$
Then the $u^{\ast}-$ parameter curves are unit speed spacelike Bertrand
geodesics and $v^{\ast}-$ parameter curves forms orthogonal spacelike
parallels. If the following properties are satisfied:

\begin{description}
\item[i)] $\left\vert \sigma-\sigma^{\ast}\right\vert =$ constant;

\item[ii)] $\sigma-\sigma^{\ast}\bot$ $b;$

\item[iii)] $\sigma-\sigma^{\ast}\bot$ $b^{\ast}$;

\item[iv)] $\left\langle b,b^{\ast}\right\rangle _{\mathbb{L}}=$ constant,
\end{description}

\noindent then $\varphi$ is called a Razzaboni transformation and the surface
$\sigma^{\ast}$ is called dual\linebreak Razzaboni surface of $\sigma.$ By the
first property, the distance between corresponding points of the surfaces is
constant. So, we can write%
\[
\sigma^{\ast}(u,v)=\varphi\left(  \sigma(u,v)\right)  =\sigma(u,v)+An(u,v).
\]
Since the geodesics of the Razzaboni surfaces $\sigma$ and $\sigma^{\ast}$ are
Bertrand mates, we have $n=n^{\ast}$. Also we get
\[
\sigma_{u}^{\ast}=\sigma_{u}+An_{u}=(1-A\kappa)t-A\tau b=\tau(Bt-Ab).
\]
That is
\[
t^{\ast}=\frac{Bt-Ab}{\sqrt{A^{2}+B^{2}}}.
\]
Since $b^{\ast}=-t^{\ast}\times_{\mathbb{L}}n^{\ast},$ we obtain%
\[
b^{\ast}=\frac{At+Bb}{\sqrt{A^{2}+B^{2}}}.
\]
And
\[
du^{\ast}=\sqrt{A^{2}+B^{2}}\tau du.
\]
Then, we have
\begin{align*}
\left\langle \sigma^{\ast}-\sigma,b\right\rangle _{\mathbb{L}} &
=\left\langle An,b\right\rangle _{\mathbb{L}}=0,\\
\left\langle \sigma^{\ast}-\sigma,b^{\ast}\right\rangle _{\mathbb{L}} &
=\frac{1}{\sqrt{A^{2}+B^{2}}}\left\langle An,At+Bb\right\rangle _{\mathbb{L}%
}=0,\\
\left\langle b,b^{\ast}\right\rangle _{\mathbb{L}} &  =\frac{1}{\sqrt
{A^{2}+B^{2}}}\left\langle b,At+Bb\right\rangle _{\mathbb{L}}=\frac{B}%
{\sqrt{A^{2}+B^{2}}}.
\end{align*}
Thus, the properties ii), iii) and iv) are also satisfied which means that
$\varphi:\sigma\rightarrow\sigma^{\ast}$ is a Razzaboni transformation. On the
other hand, we have
\[
t_{u^{\ast}}^{\ast}=\frac{B}{(A^{2}+B^{2})\tau}t_{u}-\frac{A}{(A^{2}%
+B^{2})\tau}b_{u}=\frac{-\kappa B+A\tau}{(A^{2}+B^{2})\tau}n.
\]
Then, the curvature of Bertrand geodesic of dual Razzaboni surface is found
as
\[
\kappa^{\ast}=\left\langle t_{u^{\ast}}^{\ast},n^{\ast}\right\rangle
_{\mathbb{L}}=\frac{-\kappa B+A\tau}{(A^{2}+B^{2})\tau}\left\langle
n,n\right\rangle _{\mathbb{L}}=\frac{B\kappa-A\tau}{(A^{2}+B^{2})\tau}.
\]
Moreover, we have
\[
b_{u^{\ast}}^{\ast}=\frac{A}{(A^{2}+B^{2})\tau}t_{u}+\frac{B}{(A^{2}%
+B^{2})\tau}b_{u}=\frac{-A\kappa-B\tau}{(A^{2}+B^{2})\tau}n=\frac{-1}%
{(A^{2}+B^{2})\tau}n.
\]
Then, we get
\[
\tau^{\ast}=\left\langle b_{u^{\ast}}^{\ast},n^{\ast}\right\rangle
_{\mathbb{L}}=\frac{-1}{(A^{2}+B^{2})\tau}\left\langle n,n\right\rangle
_{\mathbb{L}}=\frac{1}{(A^{2}+B^{2})\tau}.
\]
\bigskip

\noindent\textbf{Case 2: Geodesic Bertrand curves of the Razzaboni surface
have timelike binormals.}

\noindent Let $\sigma=\sigma(u,v)$ be a one parameter family of geodesic
Bertrand curve with timelike binormals in $\mathbb{M}^{3}.$ And let $u$ and
$v$ be the geodesic coordinates of the Razzaboni surface $\sigma.$ Then the
first fundamental form of the surface is of the form%
\[
I=du^{2}-\lambda^{2}dv^{2}.
\]
Here the curves $u-$ parameter curves are unit speed spacelike Bertrand
geodesics and $v-$ parameter curves forms orthogonal timelike parallels. Since
$\left\langle \sigma_{u},\sigma_{v}\right\rangle _{\mathbb{L}}=0$ and the
principal normal $n$ of geodesics are normal to the surface, then the tangent
vectors of coordinate curves are given as
\[
\sigma_{u}=t\text{ and }\sigma_{v}=\lambda b
\]
where $b$ denotes the binormal vector field of geodesics. The variation of the
Frenet Frame $\{t,n,b,\}$ of Bertrand geodesics in $u-$ direction is obtained
as follows:%
\[
\left[
\begin{array}
[c]{c}%
t\\
n\\
b
\end{array}
\right]  _{u}=\left[
\begin{array}
[c]{ccc}%
0 & \kappa & 0\\
-\kappa & 0 & \tau\\
0 & \tau & 0
\end{array}
\right]  \left[
\begin{array}
[c]{c}%
t\\
n\\
b
\end{array}
\right]
\]
by setting $(\varepsilon_{1},\varepsilon_{2},\varepsilon_{3})=(1,1,-1)$ in
equation 1. The variation of $\{t,n,b\}$ in $v-$ direction must be in the
following form:%
\[
\left[
\begin{array}
[c]{c}%
t\\
n\\
b
\end{array}
\right]  _{v}=\left[
\begin{array}
[c]{ccc}%
0 & \alpha & \beta\\
-\alpha & 0 & \gamma\\
\beta & \gamma & 0
\end{array}
\right]  \left[
\begin{array}
[c]{c}%
t\\
n\\
b
\end{array}
\right]  .
\]
The compatibility condition $\sigma_{uv}=\sigma_{vu}$ implies that%
\[
\alpha n+\beta b=(\lambda\tau)n+\lambda_{u}b.
\]
So, we get $\alpha=\lambda\tau$ and $\beta=\lambda_{u}.$ Then the variation of
$\{t,n,b\}$ in $v-$ direction becomes
\[
\left[
\begin{array}
[c]{c}%
t\\
n\\
b
\end{array}
\right]  _{v}=\left[
\begin{array}
[c]{ccc}%
0 & \lambda\tau & \lambda_{u}\\
-\lambda\tau & 0 & \gamma\\
\lambda_{u} & \gamma & 0
\end{array}
\right]  \left[
\begin{array}
[c]{c}%
t\\
n\\
b
\end{array}
\right]  .
\]
On the other hand, the compatibility condition $t_{uv}=t_{vu}$ implies that
\begin{align*}
t_{uv}  &  =(\kappa n)_{v}=(-\lambda\kappa\tau)t+(\kappa_{v})n+(\kappa
\gamma)b,\\
t_{vu}  &  =(\lambda\tau n+\lambda_{u}b)_{u}=(-\lambda\kappa\tau
)t+(\lambda\tau_{u}+2\lambda_{u}\tau)n+(\kappa\gamma)b.
\end{align*}
So we get
\[
\kappa_{v}=\lambda\tau_{u}+2\lambda_{u}\tau.
\]
And by the compatibility condition $n_{uv}=n_{vu},$ we obtain%
\begin{align*}
n_{uv}  &  =(-\kappa t+\tau b)_{v}=(-\kappa_{v}+\lambda_{u}\tau)t+(-\kappa
\lambda\tau+\tau\gamma)n+(-\lambda_{u}\kappa+\tau_{v})b,\\
n_{vu}  &  =(-\lambda\tau t+\gamma b)_{u}=(-\lambda_{u}\tau-\lambda\tau
_{u})t+(-\kappa\lambda\tau+\tau\gamma)n+(\gamma_{u})b.
\end{align*}
Thus, we get
\[
\tau_{v}=\lambda_{u}\kappa+\gamma_{u}.
\]
Finally, the compatibility condition $b_{uv}=b_{vu}$ gives that%
\begin{align*}
b_{uv}  &  =(\tau n)_{v}=(-\lambda\tau^{2})t+(\tau_{v})n+(\tau\gamma)b,\\
b_{vu}  &  =(\lambda_{u}t+\gamma n)_{u}=(\lambda_{uu}-\kappa\gamma
)t+(\lambda_{u}\kappa+\gamma_{u})n+(\tau\gamma)b.
\end{align*}
Then, we get
\[
\lambda_{uu}=-\lambda\tau^{2}+\kappa\gamma
\]
The system%
\begin{align*}
\kappa_{v}  &  =\lambda\tau_{u}+2\lambda_{u}\tau,\\
\tau_{v}  &  =\lambda_{u}\kappa+\gamma_{u},\\
\lambda_{uu}  &  =-\lambda\tau^{2}+\kappa\gamma
\end{align*}
can be considered as the Gauss-Minardi-Codazzi equations for the surface. By
the constraint
\[
A\kappa+B\tau=1,
\]
the system is well determined and the surface $\sigma$ is guaranteed to be a
Razzaboni surface.

\bigskip

\noindent Similar to case 1, we can define Razzaboni transformation
$\varphi:\sigma\rightarrow\sigma^{\ast}$ as
\[
\sigma^{\ast}(u,v)=\sigma(u,v)+An(u,v).
\]
It is easily seen that the transformation $\varphi$ satisfies all properties
of Razzaboni transformation. Here $u^{\ast}$ and $v^{\ast}$ are the geodesic
coordinates of dual Razzaboni surface $\sigma^{\ast}.$ And let $\{t^{\ast
},n^{\ast},b^{\ast}\}$ be the Frenet Frame of Bertrand geodesics of dual
Razzaboni surface $\sigma^{\ast}.$ Then the $u^{\ast}-$ parameter curves are
unit speed spacelike Bertrand geodesics and $v^{\ast}-$ parameter curves forms
orthogonal timelike parallels. Since the geodesics of the Razzaboni surface
and its dual are Bertrand mates, we have $n=n^{\ast}$. Also we have
\[
\sigma_{u}^{\ast}=\sigma_{u}+An_{u}=(1-A\kappa)t+A\tau b=\tau(Bt+Ab).
\]
That is
\[
t^{\ast}=\frac{Bt+Ab}{\sqrt{B^{2}-A^{2}}}.
\]
Since $b^{\ast}=-t^{\ast}\times_{\mathbb{L}}n^{\ast},$ we obtain%
\[
b^{\ast}=\frac{-At+Bb}{\sqrt{B^{2}-A^{2}}}.
\]
And
\[
du^{\ast}=\sqrt{B^{2}-A^{2}}\tau du.
\]
Moreover, we have
\begin{align*}
t_{u^{\ast}}^{\ast}  &  =\frac{B}{(B^{2}-A^{2})\tau}t_{u}+\frac{A}%
{(B^{2}-A^{2})\tau}b_{u}=\frac{\kappa B+A\tau}{(B^{2}-A^{2})\tau}n,\\
\kappa^{\ast}  &  =\left\langle t_{u^{\ast}}^{\ast},n^{\ast}\right\rangle
_{\mathbb{L}}=\frac{\kappa B+A\tau}{(B^{2}-A^{2})\tau}\left\langle
n,n\right\rangle _{\mathbb{L}}=\frac{\kappa B+A\tau}{(B^{2}-A^{2})\tau},\\
b_{u^{\ast}}^{\ast}  &  =-\frac{A}{(B^{2}-A^{2})\tau}t_{u}+\frac{B}%
{(B^{2}-A^{2})\tau}b_{u}=\frac{-A\kappa+B\tau}{(B^{2}-A^{2})\tau}n,\\
\tau^{\ast}  &  =\left\langle b_{u^{\ast}}^{\ast},n^{\ast}\right\rangle
_{\mathbb{L}}=\frac{-A\kappa+B\tau}{(B^{2}-A^{2})\tau}\left\langle
n,n\right\rangle _{\mathbb{L}}=\frac{-A\kappa+B\tau}{(B^{2}-A^{2})\tau}.
\end{align*}
\noindent\textbf{Case 3: Geodesic Bertrand curves of the Razzaboni surface are
timelike.}

\noindent Let $\sigma=\sigma(u,v)$ be a one parameter family of timelike
geodesic Bertrand curve in $\mathbb{M}^{3}.$ And let $u$ and $v$ be the
geodesic coordinates of the Razzaboni surface $\sigma.$ Then the first
fundamental form of the surface is of the form%
\[
I=-du^{2}+\lambda^{2}dv^{2}.
\]
Here the $u-$ parameter curves are unit speed timelike Bertrand geodesics and
$v-$ parameter curves forms orthogonal spacelike parallels. Since
$\left\langle \sigma_{u},\sigma_{v}\right\rangle _{\mathbb{L}}=0$ and the
principal normal $n$ of geodesics are normal to the surface, then the tangent
vectors of coordinate curves are given as
\[
\sigma_{u}=t\text{ and }\sigma_{v}=\lambda b
\]
where $b$ denotes the binormal vector field of geodesics. The variation of the
Frenet Frame $\{t,n,b,\}$ of Bertrand geodesics in $u-$ direction is obtained
as follows:%
\[
\left[
\begin{array}
[c]{c}%
t\\
n\\
b
\end{array}
\right]  _{u}=\left[
\begin{array}
[c]{ccc}%
0 & \kappa & 0\\
\kappa & 0 & -\tau\\
0 & \tau & 0
\end{array}
\right]  \left[
\begin{array}
[c]{c}%
t\\
n\\
b
\end{array}
\right]
\]
by setting $(\varepsilon_{1},\varepsilon_{2},\varepsilon_{3})=(-1,1,1)$ in
equation 1. The variation of $\{t,n,b\}$ in $v-$ direction must be in the
following form:%
\[
\left[
\begin{array}
[c]{c}%
t\\
n\\
b
\end{array}
\right]  _{v}=\left[
\begin{array}
[c]{ccc}%
0 & \alpha & \beta\\
\alpha & 0 & \gamma\\
\beta & -\gamma & 0
\end{array}
\right]  \left[
\begin{array}
[c]{c}%
t\\
n\\
b
\end{array}
\right]  .
\]
The compatibility condition $\sigma_{uv}=\sigma_{vu}$ implies that%
\[
\alpha n+\beta b=(\lambda\tau)n+\lambda_{u}b.
\]
So, we get $\alpha=\lambda\tau$ and $\beta=\lambda_{u}.$ Then the variation of
$\{t,n,b\}$ in $v-$ direction becomes
\[
\left[
\begin{array}
[c]{c}%
t\\
n\\
b
\end{array}
\right]  _{v}=\left[
\begin{array}
[c]{ccc}%
0 & \lambda\tau & \lambda_{u}\\
\lambda\tau & 0 & \gamma\\
\lambda_{u} & -\gamma & 0
\end{array}
\right]  \left[
\begin{array}
[c]{c}%
t\\
n\\
b
\end{array}
\right]  .
\]
On the other hand, the compatibility condition $t_{uv}=t_{vu}$ implies that
\begin{align*}
t_{uv}  &  =(\kappa n)_{v}=(\lambda\kappa\tau)t+(\kappa_{v})n+(\kappa
\gamma)b,\\
t_{vu}  &  =(\lambda\tau n+\lambda_{u}b)_{u}=(\lambda\kappa\tau)t+(\lambda
\tau_{u}+2\lambda_{u}\tau)n+(\kappa\gamma)b.
\end{align*}
So we get
\[
\kappa_{v}=\lambda\tau_{u}+2\lambda_{u}\tau.
\]
And by the compatibility condition $n_{uv}=n_{vu},$ we obtain%
\begin{align*}
n_{uv}  &  =(\kappa t-\tau b)_{v}=(\kappa_{v}-\lambda_{u}\tau)t+(\kappa
\lambda\tau+\tau\gamma)n+(\lambda_{u}\kappa-\tau_{v})b,\\
n_{vu}  &  =(\lambda\tau t+\gamma b)_{u}=(\lambda_{u}\tau+\lambda\tau
_{u})t+(\kappa\lambda\tau+\tau\gamma)n+(\gamma_{u})b.
\end{align*}
Thus, we get
\[
\tau_{v}=\lambda_{u}\kappa-\gamma_{u}.
\]
Finally, the compatibility condition $b_{uv}=b_{vu}$ gives that%
\begin{align*}
b_{uv}  &  =(\tau n)_{v}=(\lambda\tau^{2})t+(\tau_{v})n+(\tau\gamma)b,\\
b_{vu}  &  =(\lambda_{u}t-\gamma n)_{u}=(\lambda_{uu}-\kappa\gamma
)t+(\lambda_{u}\kappa-\gamma_{u})n+(\tau\gamma)b.
\end{align*}
Then, we get
\[
\lambda_{uu}=\lambda\tau^{2}+\kappa\gamma
\]
The system%
\begin{align*}
\kappa_{v}  &  =\lambda\tau_{u}+2\lambda_{u}\tau,\\
\tau_{v}  &  =\lambda_{u}\kappa-\gamma_{u},\\
\lambda_{uu}  &  =\lambda\tau^{2}+\kappa\gamma
\end{align*}
may be regarded as the Gauss-Minardi-Codazzi equations for the surface. If the
constraint
\[
A\kappa+B\tau=1
\]
is imposed then the system is well determined and the surface $\sigma$ is
guaranteed to be a Razzaboni surface.

\bigskip

\noindent Again, we define the Razzaboni transformation $\varphi
:\sigma\rightarrow\sigma^{\ast}$ as follows;%
\[
\sigma^{\ast}(u,v)=\sigma(u,v)+An(u,v).
\]
Let $u^{\ast}$ and $v^{\ast}$ be the geodesic coordinates of dual Razzaboni
surface $\sigma^{\ast}.$ And let $\{t^{\ast},n^{\ast},b^{\ast}\}$ be the
Frenet Frame of Bertrand geodesics of dual Razzaboni surface $\sigma^{\ast}.$
Then the $u^{\ast}-$ parameter curves are unit speed timelike Bertrand
geodesics and $v^{\ast}-$ parameter curves forms orthogonal spacelike
parallels. Since the geodesics of the Razzaboni surface and its dual are
Bertrand mates, we have $n=n^{\ast}$. Also we have
\[
\sigma_{u}^{\ast}=\sigma_{u}-An_{u}=(1-A\kappa)t+A\tau b=\tau(Bt+Ab).
\]
That is
\[
t^{\ast}=\frac{Bt+Ab}{\sqrt{B^{2}-A^{2}}}.
\]
Since $b^{\ast}=-t^{\ast}\times_{\mathbb{L}}n^{\ast},$ we obtain%
\[
b^{\ast}=\frac{-At+Bb}{\sqrt{B^{2}-A^{2}}}.
\]
And
\[
du^{\ast}=\sqrt{B^{2}-A^{2}}\tau du.
\]
On the other hand, we have
\begin{align*}
t_{u^{\ast}}^{\ast}  &  =\frac{B}{(B^{2}-A^{2})\tau}t_{u}+\frac{A}%
{(B^{2}-A^{2})\tau}b_{u}=\frac{B\kappa+A\tau}{(B^{2}-A^{2})\tau}n,\\
\kappa^{\ast}  &  =\left\langle t_{u^{\ast}}^{\ast},n^{\ast}\right\rangle
_{\mathbb{L}}=\frac{B\kappa+A\tau}{(B^{2}-A^{2})\tau}\left\langle
n,n\right\rangle _{\mathbb{L}}=\frac{B\kappa+A\tau}{(B^{2}-A^{2})\tau},\\
b_{u^{\ast}}^{\ast}  &  =\frac{-A}{(B^{2}-A^{2})\tau}t_{u}+\frac{B}%
{(B^{2}-A^{2})\tau}b_{u}=\frac{-A\kappa+B\tau}{(B^{2}-A^{2})\tau}n,\\
\tau^{\ast}  &  =\left\langle b_{u^{\ast}}^{\ast},n^{\ast}\right\rangle
_{\mathbb{L}}=\frac{-A\kappa+B\tau}{(B^{2}-A^{2})\tau}\left\langle
n,n\right\rangle _{\mathbb{L}}=\frac{-A\kappa+B\tau}{(B^{2}-A^{2})\tau}.
\end{align*}
\noindent\textbf{Conclusion}

\noindent In first case, we obtain that the solutions $\kappa,$ $\tau,$
$\lambda$ and $\gamma$ of the system
\begin{align*}
\kappa_{v}  &  =\lambda\tau_{u}+2\lambda_{u}\tau,\\
\tau_{v}  &  =-\lambda_{u}\kappa-\gamma_{u},\\
\lambda_{uu}  &  =-\lambda\tau^{2}-\kappa\gamma,\\
A\kappa+B\tau &  =1
\end{align*}
constitutes the spacelike Razzaboni surfaces in $\mathbb{M}^{3}.$ In the
second case we obtain that the solutions $\kappa,$ $\tau,$ $\lambda$ and
$\gamma$ of the system
\begin{gather*}
\kappa_{v}=\lambda\tau_{u}+2\lambda_{u}\tau,\\
\tau_{v}=\lambda_{u}\kappa+\gamma_{u},\\
\lambda_{uu}=-\lambda\tau^{2}+\kappa\gamma\\
A\kappa+B\tau=1
\end{gather*}
constitutes the timelike Razzaboni surfaces with timelike $v-$ parameter
curves in $\mathbb{M}^{3}.$

\bigskip

\noindent For first and second case, the second fundamental form of the
surface $\ \sigma$ is of the form
\[
II=-\kappa du^{2}-2\lambda\tau dudv+\frac{\lambda}{\kappa}(-\lambda
_{uu}-\lambda\tau^{2})dv^{2}.
\]
Then the Gaussian curvature of the surface is obtained as
\[
K=-\frac{\lambda_{uu}}{\lambda}.
\]

\bigskip

\noindent In the last case, we obtain that the solutions $\kappa,$ $\tau,$
$\lambda$ and $\gamma$ of the system%
\begin{gather*}
\kappa_{v}=\lambda\tau_{u}+2\lambda_{u}\tau,\\
\tau_{v}=\lambda_{u}\kappa-\gamma_{u},\\
\lambda_{uu}=\lambda\tau^{2}+\kappa\gamma\\
A\kappa+B\tau=1
\end{gather*}
constitutes the timelike Razzaboni surfaces with timelike $u-$ parameter
curves in $\mathbb{M}^{3}.$ In this case, the second fundamental form of the
surface is of the form
\[
II=-\kappa du^{2}-2\lambda\tau dudv+\frac{\lambda}{\kappa}(-\lambda
_{uu}+\lambda\tau^{2})dv^{2}.
\]
Then the Gaussian curvature of the surface is obtained as
\[
K=\frac{\lambda_{uu}}{\lambda}.
\]

\bigskip

\noindent In the case of $A=0,$ which means that Bertrand geodesics of the
surface $\sigma$ must have constant torsion
\[
\tau=\frac{1}{B},
\]
then Razzaboni transformation of the surface coincides the main Razzaboni
surface. That is
\[
\varphi(\sigma)=\sigma.
\]

\bigskip

\noindent In case of $B=0,$ which means that Bertrand geodesics of the surface
$\sigma$ must have constant curvature, we have
\[
\kappa=\frac{1}{A}.
\]
Thus, Razzaboni transformation is defined as
\[
\varphi(\sigma(u,v))=\sigma^{\ast}(u,v)=\sigma(u,v)+\frac{1}{\kappa}n(u,v).
\]
Since
\[
\kappa^{\ast}=-\frac{1}{A}=-\kappa,
\]
then, Razzaboni transformation maps the surface $\sigma,$ whose Bertrand
geodesics have constant curvature, to the surface $\sigma^{\ast}$ whose
Bertrand geodesics also have constant curvature with opposite sign. Moreover,
the torsion of Bertrand geodesics of the surfaces satisfy the relation
\[
\tau\tau^{\ast}=\frac{1}{A^{2}}.
\]

\end{document}